\documentclass[12pt,a4paper]{article}
\usepackage{amsmath,amsfonts,amssymb}
\newtheorem{theo}{Theorem}[section]
\newtheorem{prop}[theo]{Proposition}

\newtheorem{lem}[theo]{Lemma}

\def\genschur{\mathfrak{G}}
\def\X{\mathbb X}

\def\A{A}
\def\B{B}
\def\N{\mathbb N}
\def\Z{\mathbb Z}
\def\C{\mathcal C}

\def\V{\widetilde V}
\def\l{l}

\def\pf{\noindent {\it Proof.\ }}
\def\qed{\hfill \rule{4pt}{7pt}}
\def\bsquare{\rule{5pt}{5pt}}
\def\ybox{\framebox(0.6,0.6){}}
\def\xbox{\framebox(0.6,0.6){$\star$}}

\begin{document}
\begin{center}
{\large \bf Division and the Giambelli Identity}
\end{center}

\begin{center}
{\small
Susan Y. J. Wu$^1$ and Arthur L. B. Yang$^2$\\[5pt]
Center for Combinatorics, LPMC\\
Nankai University, Tianjin 300071, P. R. China\\
E-mail: $^1$wuyijun@126.com and $^2$yang@nankai.edu.cn}
\end{center}

\vspace{0.3cm} \noindent{\bf Abstract.} Given two polynomials
$f(x)$ and $g(x)$, we extend the formula expressing the remainder
in terms of the roots of these two polynomials to the case where
$f(x)$ is a Laurent polynomial. This allows us to give new
expressions of a Schur function, which generalize the Giambelli
identity.



\noindent {\bf Corresponding Author:} Arthur L. B. Yang,
yang@nankai.edu.cn

\section{Introduction}

The Euclidean algorithm is an algorithm to determine the greatest
common divisor of two integers, which appeared in \emph{Euclid's
Elements} around 300 BC. However it is easily generalized to
polynomials in one variable $x$ over the field of real numbers. It
turns out that this process generates symmetric functions over the
variable sets $A$ and $B$, if $A$ and $B$ are the alphabets of
roots of the two original polynomials. By developing this point of
view in \cite{cbms}, Lascoux obtained the explicit expressions of
remainders in terms of Schur functions.

We assume that the reader is familiar with the background of the
theory of symmetric functions \cite{cbms, Macd, Stanley}. We use
nondecreasing partitions to index Schur functions. Let $\A$ be of
cardinality $n,\, I\in\N^n$ be a partition contained in some
rectangular partition $\square=m^n$, and $J$ be the complementary
partition of $I$ in $\square$. We denote the set $\{a^{-1}:
a\in\A\}$ by $\A^{\vee}$. Let $u=a_1\cdots a_n$ be the product of
all the variables in $\A$. Taking the expression of a Schur
function in terms of the Vandermonde matrix (\cite[p. 40]{Macd}),
then one has the following relation between the Schur functions in
$\A$ and those in $\A^{\vee}$:
\begin{equation}\label{eq-inv-self}
S_{I}(\A^{\vee})={S_{J}(\A)}u^{-m}.
\end{equation}

Taking an extra indeterminate $z$ and two alphabets $\A, \B$, then
the \emph{complete symmetric functions} $S^k(\A-\B)$ are defined
by the generating function
\begin{equation}
\sum_{k\geqslant 0}S^k(\A-\B)z^k={\frac{\prod_{b\in\B}(1-bz)}{
\prod_{a\in\A}(1-az)}}.
\end{equation}
Given two sets of alphabets $\{ \A_1, \A_2, \ldots, \A_n \}$ and
$\{ \B_1, \B_2, \ldots, \B_n\}$, and $I,J\in \N^n$,  then the
\emph{multi-Schur function} of index $J/I$ is defined as follows
\cite{cbms}:
\begin{equation}
S_{J/I}({\A}_1-{\B}_1; \ldots; {\A}_n-{\B}_n)
    := \bigl| \, S_{j_k -i_{\l}+k-\l}({\A}_k -{\B}_k)
           \, \bigr|_{1\leq \l,k\leq n}.
\end{equation}
If each column has the same argument $\A-\B$, we denote the
multi-Schur function by $S_{J/I}(\A-\B)$.

Lascoux \cite{cbms} proved that
\begin{theo}\label{l-thm}
The $r$-th remainder in the division of $S^m(x-B)$ by $S^n({x-A})$
is equal to
\begin{equation}
S_{1^{n-r};(m-n+r)^r}(A-x,A-B).
\end{equation}
\end{theo}

In section 2, we adapt division to the case of the division of a
Laurent polynomial by a usual polynomial, and we give several
expressions of the first remainder as a Schur function. The
Lagrange interpolation and Lagrange functional are used to
reconstruct these remainders. To proceed the Euclidean algorithm,
Theorem \ref{l-thm} allows us to obtain expressions for other
remainders in terms of Schur functions.

For an arbitrary Schur function of shape $J$,  the Giambelli
identity provides a formula which expresses $S_{J}(A)$ as a
determinant with entries being Schur functions of hook shapes
\cite{Gi, Macd}. Many combinatorial proofs and extensions of the
Giambelli identity have appeared, and we refer the reader to
\cite{ER, FK, St}. By expressing the remainders of $x^{k},\,k\in
\N$ by $S^n({x-A})$ as Schur functions, Lascoux presents another
proof for the Giambelli identity \cite{cbms}. We find that this
idea can also be used to study the extension of Schur functions
with negative indices \cite{HouMu}, denoted
$\genschur_J(\A),\,J\in \Z^n$, which are needed when interpreting
them as characters of the linear group. Following the treatment of
Lascoux in Section 3, we construct a matrix with column indices in
$\Z$, that we call double companion matrix, by putting the
coefficients of the remainders of all $x^{k}$, $k\in \Z$ into this
matrix. Taking minors of this matrix, we obtain new determinantal
expressions for $\genschur_J(\A)$, which generalize the usual
Giambelli identity. We should point out that this extension of the
Giambelli identity can also be derived from the following theorem
given by Hou and Mu \cite{HouMu}
\begin{theo}
 Given $n$ recurrent sequences
$T^{(i)}=\{T^{(i)}_m:m\in \Z\}\, (1\leq i\leq n)$ with the same
characteristic polynomial having the root set $A$, then we have
\begin{equation}\label{det-rec}
\genschur_J(\A)=\frac{|T^{(k)}_{j_l+l-1}|_{1\leq k,l\leq
n}}{|T^{(k)}_{l-1}|_{1\leq k,l\leq n}}.
\end{equation}
\end{theo}

\section{Division}\label{sec-1}

Given two polynomials $f(x)$ and $g(x)$, there exists a unique
pair $(q(x), r(x))$ such that
\begin{equation}\label{eq-div-def}
f(x)=q(x)g(x)+r(x)\quad and \quad \deg(r(x))<\deg(g(x)),
\end{equation}
where we denote the degree of a polynomial by $\deg(\,)$.

Equation \eqref{eq-div-def} remains valid if $f(x)$ and $q(x)$ are
polynomials in $x^{-1}$, i.e. there exists a unique polynomial
$r(x)$ of degree $<n$, that we still call the remainder.

In the case of a general Laurent polynomial, one would uniquely
decompose it into $f_1(x)+f_2(x^{-1})$, with $f_2(0)=0$. Formulas
for the remainders in the case of polynomials are well known, and
we shall show how to adapt them to the case where $f(x)$ is a
polynomial in $x^{-1}$.

Given two sets of variables $\A$ and $\B$, denote by $R(\A,\B)$
the product $\prod_{a\in{\A},b\in{\B}} (a-b)$, and by ${\A}-{\B}$
the set difference. Supposing $g(x)$ to be monic, with set of
roots $\A=\{a_1,\, a_2,\, \ldots,\, a_n\}$ (that we suppose
distinct), then we can write it $g(x)=R(x, \A)$. In terms of $\A$,
the remainder $r(x)$ is characterized by the conditions
\begin{equation}\label{condition}
\left\{
\begin{array}{lcl}
r(a)& = & f(a), \quad \mbox{for each $a\in\A$}\\
\deg(r(x)) & \leq & n-1.
\end{array}
\right.
\end{equation}

A polynomial of degree less than $n$ is determined by its values
in $n$ points. One can reconstruct it by the Lagrange formula,
that we shall interpret with the help of a \emph{Lagrange
functional} $L_{\A}$ \cite{cbms}. Let $\mathfrak{Sym}(\A)$ be the
ring of symmetric functions in $\A$, and let
$\mathfrak{Sym}(1|n-1)$ be the space of Laurent polynomials of a
set $\X$ of  $n$ variables $\{x_1,\,x_2\,\ldots,\,x_n\}$, which
are symmetrical in the last $n-1$ variables. Then $L_{\A}$ is
defined by
\begin{equation}
\mathfrak{Sym}(1|n-1)\ni p \longrightarrow
L_{\A}(p):=\sum_{a\in\A}{\frac{p(a,\A-a)}{
R(a,\A-a)}}\in\mathfrak{Sym}(\A).
\end{equation}
In terms of $L_{\A}$, the expression of the remainder is
\begin{equation}\label{eq-lag-expr}
r(x)=L_{\A}\biggl(r(x_1)R(x,\X-x_1)\biggr).
\end{equation}

The main theorem is
\begin{theo}\label{main-thm}
Given $k\in\N$ and $\A$ of cardinality $n$, then the remainder of
$x^{-k}$ modulo by $R(x, \A)$ is equal to

\rm{(i)} \quad ${S_{k^{n-1}}(\A-x)}{u^{-k}}$;

\rm{(ii)} \quad $(-1)^{n-1}x^{n-1}S_{1^{n-1};
k}(\A^{\vee}-x^{-1};\A^{\vee})$;

\rm{(iii)} \quad Given $\B$ of cardinality $m$, the remainder of
$R(x^{-1}, \B)$ is equal to\\
\vskip -5mm \quad\quad\quad $(-1)^{n-1}x^{n-1}S_{1^{n-1};
m}(\A^{\vee}-x^{-1};\A^{\vee}-\B)$.
\end{theo}
\pf (i) \quad The polynomial ${S_{k^{n-1}}(\A-x)}$ is of degree
$\leq n-1$ because $x$ appears in degree $1$ in each column.
Specializing it into any element of $\A$, say $x=a_1$, we get
${S_{k^{n-1}}(\A-x)}{u^{-k}}=(a_2\cdots a_n)^ku^{-k}=a_1^{-k}$,
and therefore this polynomial is the remainder of $x^{-k}$.

(ii) \quad We expand the Schur function by linearity on $x^{-1}$,
and obtain
$$
\begin{array}{rcl}
(-1)^{n-1}x^{n-1}S_{1^{n-1}; k}(\A^{\vee}-x^{-1};\A^{\vee}) & = &
\sum_{\l=0}^{n-1}(-1)^{n-1+\l}x^{n-1-\l} S_{1^{n-1-\l},
k}(\A^{\vee})\\[5pt]
&=& \sum_{\l=0}^{n-1}(-x)^{\l} S_{1^\l, k}(\A^{\vee})\\[5pt]
&=& \sum_{\l=0}^{n-1}(-x)^{\l} S_{(k-1)^\l, k^{(n-1)-\l}}(\A)u^{-k}\\[5pt]
&=& {S_{k^{n-1}}(\A-x)}u^{-k},
\end{array}
$$
the third step using \eqref{eq-inv-self}.

(iii)\quad By linearity (ii) implies (iii), but let us check it
directly using the Lagrange interpolation. Thanks to
\eqref{condition} and \eqref{eq-lag-expr}, we have
\begin{equation}\label{eq-1}
r(x)=L_{\A}\biggl(R(x_1^{-1},\B)R(x,\X-x_1)\biggr).
\end{equation}
Let $\Delta(\A)=\prod_{1\leq i<j\leq n}(a_j-a_i)$. Since for any
$k\in\N$,
\begin{eqnarray*}
L_{\A}(x_1^{-k})&=& \sum_{a\in\A}{\frac{a^{-k}}{
R(a,\A-a)}}\\[5pt]
&=&{\frac{1}{\Delta(\A)}}
\begin{vmatrix}
a_1^0&a_1^1&\cdots &a_1^{n-2}&a_1^{-k}\\
\vdots&\vdots& &\vdots&\vdots\\
a_n^0&a_n^1&\cdots &a_n^{n-2}&a_n^{-k}\\
\end{vmatrix}\\[5pt]
&=&(-1)^{n-1}{\frac{u^{-k}}{\Delta(\A)}}
\begin{vmatrix}
a_1^{0}&a_1^k&a_1^{k+1}&\cdots &a_1^{k+n-2}\\
\vdots&\vdots&\vdots& &\vdots\\
a_n^{0}&a_n^k&a_n^{k+1}&\cdots &a_n^{k+n-2}\\
\end{vmatrix}\\[5pt]
&=&(-1)^{n-1}u^{-k}S_{(k-1)^{n-1}}(\A)\\[5pt]
&=&(-1)^{n-1}u^{-1}S_{k-1}(\A^{\vee}),
\end{eqnarray*}
then
\begin{equation}\label{eq-2}
L_{\A}\biggl(S_k(x_1^{-1}-\B)\biggr)
=(-1)^{n-1}u^{-1}S_{k-1}(\A^{\vee}-\B).
\end{equation}
Moreover we have
\begin{equation*}
R(x_1^{-1},\B)R(x,\X-x_1) =(-1)^{n-1}{\frac{x^{n-1}}{
x_1^{-1}\cdots
x_n^{-1}}}S_{m+1}(x_1^{-1}-\B)S_{n-1}(x^{-1}-\X^{\vee}+x_1^{-1}),
\end{equation*}
which is equal to
\begin{equation}\label{eq-3}
R(x_1^{-1},\B)R(x,\X-x_1)={\frac{x^{n-1}}{ x_1^{-1}\cdots
x_n^{-1}}}S_{1^{n-1};m+1}(\X^{\vee}-x^{-1},x_1^{-1}-\B).
\end{equation}
Thus the equations \eqref{eq-1}, \eqref{eq-2} and \eqref{eq-3}
lead to
\begin{equation*}
r(x)=(-1)^{n-1}x^{n-1}S_{1^{n-1}; m}(\A^{\vee}-x^{-1};\A^{\vee}-B)
\end{equation*}
\qed

\section{The Giambelli identity}\label{sec-2}

We modify the definition of a Schur function (see also Hou and Mu
\cite{HouMu}), and for $J\in \mathbb{Z}^{n}$ put
\begin{equation}\label{gen-def}
\genschur_J(\A)={\frac{\bigl|a_k^{j_\l+\l-1}\bigr|_{1\leq \l,k\leq
n}}{\bigl|a_k^{\l-1}\bigr|_{1\leq \l,k\leq n}}}.
\end{equation}
In the case where $J\in\N^n$, it coincides with the usual
definition of the Schur function $S_J(\A)$. However, when $\A$ has
two letters, the usual Schur function $S_{4,-2}(A)$, defined as a
determinant of complete functions, is null, but
$\genschur_{4,-2}(\A)$ is not. In fact, one can get rid of
negative powers by multiplication by $u=a_1\cdots a_n$, then
$\genschur_J(\A)$ can be written as a Schur function in $\A$, as
well as in $\A^{\vee}$, up to powers of $u$. The following
property is easy to check:
\begin{lem}
For any $J\in \N^n$,
\begin{equation}
\genschur_J(\A)=S_J(\A)\quad and \quad
\genschur_{-J}(\A)=S_{J^\omega}(\A^{\vee}),
\end{equation}
where
$$-J=(-j_1, \ldots, -j_n)\quad and \quad J^{\omega}=(j_n, \ldots, j_1).$$
\end{lem}

The usual companion matrix, finite or infinite, is the matrix of
coefficients of the remainders of $x^1,\,\ldots,\,x^n$ (resp.
$x^0,\,x^1,\,\ldots,\,x^{\infty}$). We define the \emph{double
companion matrix} $\mathcal{C}(\A)$ to be the matrix of
coefficients of the remainders of
$\ldots,\,x^{-2},\,x^{-1},\,x^0,\,x^1,\,\ldots$ in the basis
$x^0,\,x^1,\,\ldots,\,x^{n-1}$, modulo $R(x,\A)$. Explicitly, for
any $k\in \Z$, if the remainder $r(x)$ of $x^k$ modulo $R(x, \A)$
is
\begin{equation}
r(x)=c_{0,k}x^0+c_{1,k}x^1+\cdots+c_{n-1,k}x^{n-1}.
\end{equation}
then we let
\begin{equation}
\mathcal{C}(\A)=\left(c_{\l-1,k}\right)_{1\leq\l\leq n, k\in Z}.
\end{equation}

For $k\in \N$, the remainder $r(x)$ of $x^k$ modulo $R(x, \A)$ is
given in \cite{cbms}
\begin{equation}
r(x)=(-1)^{n-1}S_{1^{n-1}; k-n+1}(\A-x, \A).
\end{equation}
Expanding the first $n-1$ columns according to
$S_{j}(\A-x)=S_j(\A)-xS_{j-1}(\A)$, we get
\begin{equation}
r(x)=\sum_{\l=1}^{n}(-1)^{n-\l}x^{\l-1} S_{1^{n-\l}, k-n+1}(\A).
\end{equation}
Thus for any $\l: 1\leq \l\leq n$ and $k\in \N$, we have
\begin{eqnarray}\label{r-1}
{c_{\l-1,k}}&=&(-1)^{n-\l}S_{1^{n-\l},
k-n+1}(\A)\nonumber\\[5pt]
&=&S_{k-\l+1,0^{n-\l}}(\A)=S_{0^{\l-1},k-\l+1,0^{n-\l}}(\A)=\genschur_{0^{\l-1},k-\l+1,0^{n-\l}}(\A).
\end{eqnarray}
By Theorem \ref{main-thm} the remainder $r(x)$ of $x^{-k}$ modulo
$R(x, \A)$ is
\begin{equation}
r(x)=(-1)^{n-1}x^{n-1}S_{1^{n-1}; k}(\A^{\vee}-x^{-1};\A^{\vee}).
\end{equation}
Expanding the above Schur function, we get
\begin{equation}
r(x)=\sum_{\l=1}^n (-1)^{\l-1}x^{\l-1}S_{1^{\l-1},k}(\A^{\vee}).
\end{equation}
Therefore for any $\l: 1\leq \l\leq n$ and $k\in \N$,
\begin{eqnarray}\label{r-2}
{c_{\l-1,-k}}&=&(-1)^{\l-1}S_{1^{\l-1},k}(\A^{\vee})\nonumber\\[5pt]
&=&S_{0^{n-\l},k+\l-1,0^{\l-1}}(\A^{\vee})=\genschur_{0^{\l-1},-k-\l+1,0^{n-\l}}(\A).
\end{eqnarray}

Combining equation \eqref{r-1} and \eqref{r-2}, we get
\begin{equation}
\mathcal{C}(\A)=\biggl(\genschur_{0^{\l-1},k-\l+1,0^{n-\l}}(\A)\biggr)_{1\leq\l\leq
n, k\in Z}.
\end{equation}

For any $I=[i_1,\, i_2,\, \ldots,\, i_n]\in \N^n$, let
$\mathcal{C}_I(\A)$ be the submatrix of $\mathcal{C}(\A)$ on
columns $i_1 +0, i_2+1, \ldots , i_n+n-1$. The usual companion
matrix is $\mathcal{C}_{1^n}(\A)$. The following proposition is
implicit in \cite{HouMu}.
\begin{prop} For any $m\in \Z$,
\begin{equation}
\left(\mathcal{C}_{1^n}(\A)\right)^{m}=\mathcal{C}_{m^n}(\A).
\end{equation}
\end{prop}

One can similarly define the {\it double Vandermonde matrix}:
\begin{equation*}
\V(\A):=\begin{bmatrix}
\cdots&a_1^{-2}&a_1^{-1}&a_1^0&a_1^1&a_1^2&\cdots\\
\cdots&a_2^{-2}&a_2^{-1}&a_2^0&a_2^1&a_2^2&\cdots\\
\cdots&\vdots&\vdots&\vdots&\vdots&\vdots&\cdots\\
\cdots&a_n^{-2}&a_n^{-1}&a_n^0&a_n^1&a_n^2&\cdots\\
\end{bmatrix}.
\end{equation*}
The usual Vandermonde matrix $V_{0}(\A)$  of order $n$ is the
submatrix of $\V(\A)$ on columns $0,\,1,\,\ldots,\,n-1$.

\begin{prop}\label{von-comp} Let $V_{0}(\A)$ be the finite
Vandermonde matrix on $\A$. Then
\begin{equation}
V_{0}(\A)\C(\A)=\V(\A)
\end{equation}
\end{prop}

This factorization implies that for any $J$,
$|\mathcal{C}_J(\A)V_{0}(\A)|$ is equal to the minor of $\V(\A)$
on columns $j_1 +0, j_2+1, \ldots , j_n+n-1$. Thanks to
\eqref{gen-def}, we therefore obtain the following theorem, which
generalizes Giambelli's identity to the Schur function
$\genschur_J(\A)$ (see \cite{Gi} and \cite[p. 47]{Macd}).
\begin{theo}\label{main-thm2}
\begin{equation}\label{main-eq}
\genschur_J(\A)=\bigl| \genschur_{0^{\l-1}, j_k+k-\l,0^{n-\l}}(\A)
\bigr|_{1\leq \l,k\leq n}.
\end{equation}
\end{theo}
This theorem follows also from \cite[Theorem 4.4]{HouMu} once we
check that for each $\l: 1\leq \l\leq n$,\,
$\{\genschur_{0^{\l-1}, k-\l+1,0^{n-\l}},k\in \Z\}$ is a recurrent
sequence with characteristic polynomial $R(x, \A)$.

For any weakly increasing sequence $J\in \Z^n$, let
$J_1=(j_1,\ldots,j_t)$ be the negative part and
$J_2=(j_{t+1},\ldots,j_n)$ nonnegative part. Let $(\alpha|\beta)$
be the Frobenius decomposition into diagonal hooks of
$-J_1^{\omega}$ (with rank $r_1$), and let $(\gamma|\delta)$ be
the  Frobenius decomposition of $J_2$ (with rank $r_2$) \cite[p.
3]{Macd}. Let $i\, \&\, j$ denote the partition $(1^{j},i+1)$ for
$i,j\in \N$.

Some modification on the determinant in \eqref{main-eq}
(suppressing columns having only one occurrence of $1$, the other
entries being $0$) leads to the following combinatorial version of
Theorem \ref{main-thm2}
\begin{theo} For any weakly increasing sequence $J\in\Z^n$, let $\alpha, \beta, \gamma,
\delta$ be defined as above, then
\begin{equation}\label{com}
\genschur_J(\A) =\begin{vmatrix} P&Q\\M&N
\end{vmatrix},
\end{equation}
where
$$\begin{array}{ll}
P=\left(S_{\alpha_{r_1+1-j}\&\beta_{r_1+1-i}}(\A^{\vee})\right)_{r_1\times
r_1},
&Q=\left(S_{\gamma_{j}\&(n-1-\beta_{r_1+1-i})}(\A)\right)_{r_1\times
r_2},\\[5pt]
M=\left(S_{\alpha_{r_1+1-j}\&(n-1-\delta_{i})}(\A^{\vee})\right)_{r_2\times
r_1},
&N=\left(S_{\gamma_{j}\&\delta_{i}}(\A)\right)_{r_2\times r_2}.\\
\end{array}$$
\end{theo}

For example, for $n=6$, $J=[-4,-3,-2,1,3,4]$, one has
$$\genschur_J(\A)=\begin{vmatrix}
S_{12}(\A^{\vee})    &S_{14}(\A^{\vee})    & S_{1^4,4}(A) & S_{1^4,2}(A)\\
S_{112}(\A^{\vee})   &S_{114}(\A^{\vee})   & S_{1^3,4}(A) & S_{1^3,2}(A)\\
S_{1^3,2}(\A^{\vee}) &S_{1^3,4}(\A^{\vee}) & S_{114}(A)   & S_{112}(A)\\
S_{1^5,2}(\A^{\vee}) &S_{1^5,4}(\A^{\vee}) & S_{4}(A)     & S_{2}(A)\\
\end{vmatrix}.
$$
Notice that the first two columns involve $\A^{\vee}$, and the
last two columns involve $\A$.
\begin{figure}[h,t]
\begin{center}
\begin{picture}(140,160)
\setlength{\unitlength}{8pt}

\multiput(-14,12)(1,0){3}{\ybox}
\multiput(-13,11)(1,0){1}{\ybox}
\put(-11,12){\xbox}
\put(-12,11){\xbox}
\multiput(-11,11)(0,-1){2}{$\bsquare$}
\multiput(-12,10)(0,-1){1}{$\bsquare$}
\multiput(-10,9)(0,-1){2}{$\bsquare$}
\multiput(-9,8)(1,0){2}{\ybox}
\multiput(-10,7)(1,0){4}{\ybox}
\put(-9,8){\xbox}
\put(-10,7){\xbox}


\multiput(-3.5,19)(1,0){4}{\ybox}
\multiput(-0.5,17)(0,1){2}{$\bsquare$}
\multiput(-2.5,14)(1,0){2}{\ybox}
\multiput(-1.5,13)(1,0){1}{$\bsquare$} \put(-0.5,19){\xbox}
\put(-1.5,14){\xbox}

\put(1,18){$\rightarrow$} \put(1,13.5){$\rightarrow$}
\put(1,6){$\rightarrow$} \put(1,1){$\rightarrow$}

\multiput(5,19)(0,-1){3}{$\bsquare$}
\multiput(5,15)(0,-1){4}{$\bsquare$}

\multiput(-3.5,7)(0,-1){2}{$\bsquare$}
\multiput(-3.5,5)(1,0){4}{\ybox} \multiput(-2.5,1)(1,0){2}{\ybox}
\put(-3.5,5){\xbox} \put(-2.5,1){\xbox}

\multiput(5,5)(0,1){3}{$\bsquare$}
\multiput(5,-1)(0,1){5}{$\bsquare$}

\put(8,9.8){$\Rightarrow$}

\multiput(11,18)(1,0){2}{\ybox}
\multiput(11,14)(1,0){2}{\ybox}
\multiput(11,9)(1,0){2}{\ybox}
\multiput(11,4)(1,0){2}{\ybox}
\multiput(12,17)(0,1){1}{$\bsquare$}
\multiput(12,12)(0,1){2}{$\bsquare$}
\multiput(12,6)(0,1){3}{$\bsquare$}
\multiput(12,-1)(0,1){5}{$\bsquare$}
\put(12,18){\xbox}
\put(12,14){\xbox}
\put(12,9){\xbox}
\put(12,4){\xbox}

\multiput(14,18)(1,0){4}{\ybox}
\multiput(14,14)(1,0){4}{\ybox}
\multiput(14,9)(1,0){4}{\ybox}
\multiput(14,4)(1,0){4}{\ybox}
\multiput(17,17)(0,1){1}{$\bsquare$}
\multiput(17,12)(0,1){2}{$\bsquare$}
\multiput(17,6)(0,1){3}{$\bsquare$}
\multiput(17,-1)(0,1){5}{$\bsquare$}
\put(17,18){\xbox}
\put(17,14){\xbox}
\put(17,9){\xbox}
\put(17,4){\xbox}

\multiput(22,16)(1,0){4}{\ybox}
\multiput(22,11)(1,0){4}{\ybox}
\multiput(22,5)(1,0){4}{\ybox}
\multiput(22,1)(1,0){4}{\ybox}
\multiput(22,17)(0,1){4}{$\bsquare$}
\multiput(22,12)(0,1){3}{$\bsquare$}
\multiput(22,6)(0,1){2}{$\bsquare$}
\put(22,16){\xbox}
\put(22,11){\xbox}
\put(22,5){\xbox}
\put(22,1){\xbox}

\multiput(27,16)(1,0){2}{\ybox}
\multiput(27,11)(1,0){2}{\ybox}
\multiput(27,5)(1,0){2}{\ybox}
\multiput(27,1)(1,0){2}{\ybox}
\multiput(27,17)(0,1){4}{$\bsquare$}
\multiput(27,12)(0,1){3}{$\bsquare$}
\multiput(27,6)(0,1){2}{$\bsquare$}
\put(27,16){\xbox}
\put(27,11){\xbox}
\put(27,5){\xbox}
\put(27,1){\xbox}

\put(10,-1.5 ){\line(0,1){23}}
\put(29.5,-1.5 ){\line(0,1){23}}

\thinlines \multiput(4.1,8)(1,0){3}{\line(1,0){0.5}}
\multiput(4.1,8)(0,-1){10}{\line(0,-1){0.5}}
\multiput(6.58,8)(0,-1){10}{\line(0,-1){0.5}}
\multiput(4.1,-1.5)(1,0){3}{\line(1,0){0.5}}

\multiput(4.1,20.1)(1,0){3}{\line(1,0){0.5}}
\multiput(4.1,20.1)(0,-1){9}{\line(0,-1){0.5}}
\multiput(6.58,20.1)(0,-1){9}{\line(0,-1){0.5}}
\multiput(4.1,11.6)(1,0){3}{\line(1,0){0.5}}

\multiput(10.5,10.3)(1,0){19}{\line(1,0){0.5}}
\multiput(20,-1)(0,1){22}{\line(0,1){0.5}}

\setlength{\unitlength}{1pt} \linethickness{0.5pt}
\put(-45,85){\line(0,1){65}} \put(-45,65){\line(0,-1){65}}
\put(-40,150){\oval(10,20)[tl]} \put(-40,0){\oval(10,20)[bl]}
\put(-50,85){\oval(10,20)[br]} \put(-50,65){\oval(10,20)[tr]}

\end{picture}
\end{center}
\caption{Combinatorial visualization of generalized Giambelli
identity}\label{figure-1}
\end{figure}
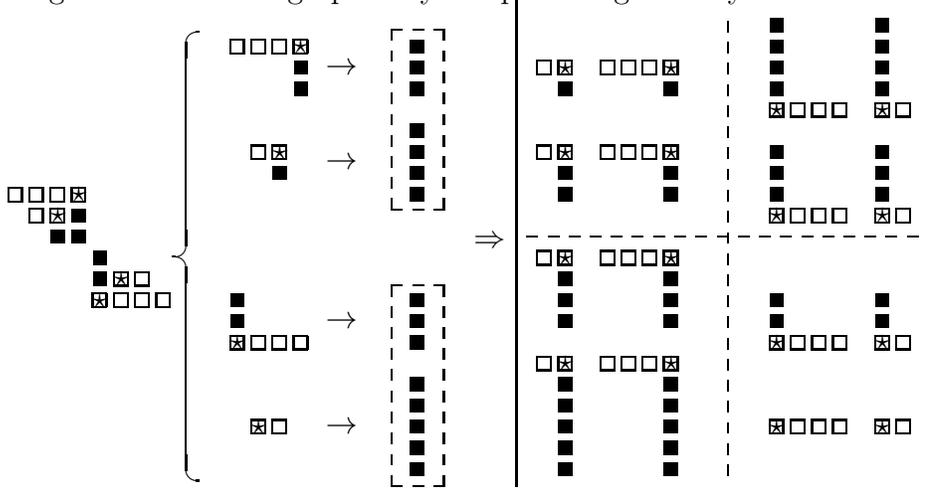
Figure \ref{figure-1} illustrates graphically the preceding
identity.

The Giambelli identity of Schur functions has been generalized in
many different ways. Lascoux and Pragacz \cite{LP} express Schur
functions as determinants of ribbon Schur functions. Hamel and
Goulden \cite{HG} use planar decompositions of skew shape tableaux
into strips, to which they associate determinantal expressions of
skew Schur functions.

Notice that in the two diagonal blocks, we have the usual
Giambelli determinants for $S_{234}(\A^{\vee})$ and $S_{134}(\A)$,
but the two other blocks are not $\mathbf{0}$, because our
function is not $S_{444/12}(\A)S_{134}(\A)$.

\vspace{.1cm} \noindent{\bf Acknowledgments.} This work was done
under the auspices of the 973 Project on Mathematical
Mechanization, the Ministry of Education, the Ministry of Science
and Technology, and the National Science Foundation of China. We
thank Professor Alain Lascoux for his useful comments, and we also
thank Dr. Q.-H. Hou and Y.-P. Mu for their helpful discussion.

\end{document}